\documentclass[11pt,letterpaper]{article}
\usepackage{times}

\pdfpagewidth=\paperwidth
\pdfpageheight=\paperheight

\usepackage{wrapfig}
\usepackage{amsmath,amssymb,latexsym,amsfonts}
\usepackage{bm}
\usepackage{stmaryrd}
\usepackage[sf,bf,compact,topmarks,calcwidth,pagestyles]{titlesec}
\usepackage[pdftex,usenames,dvipsnames]{color}
\usepackage{enumerate}
\usepackage{graphics,graphicx}
\usepackage{mathtools}
\usepackage{dialogue}
\usepackage{paralist}
\usepackage{enumitem}
\usepackage{acronym}
\usepackage{tikz}
\usepackage{pgfgantt}
\usepackage{sidecap}
\usepackage{threeparttable}
\usepackage[greek,english]{babel}
\usepackage{cite}

\usepackage{tabu}

\usepackage[colorinlistoftodos]{todonotes}
\usepackage{comment}
\usepackage{soul,color}

\usepackage{bbm}
\usepackage{mathrsfs}
\usepackage{verbatim}
\usepackage{amsmath}
\usepackage{amsfonts}
\usepackage{bm}
\usepackage{graphicx}
\usepackage{caption}
\usepackage{float}

\usepackage{amsthm}
\usepackage{xcolor}
\usepackage[colorinlistoftodos]{todonotes}
\usepackage{mathtools}
\usepackage{cite}
\usepackage{balance}

\theoremstyle{definition}

\newtheorem{problem}{Problem}

\theoremstyle{plain}

\let\oldmarginpar\marginpar
\renewcommand\marginpar[1]{\-\oldmarginpar[\raggedleft\footnotesize #1]%
	{\raggedright\footnotesize #1}}

\date{}

\usepackage{fancyhdr}
\renewcommand{\headrulewidth}{1.5pt}

\usepackage[pdftex,breaklinks=true,bookmarks=true]{hyperref}

\acrodef{cps}[\textsc{cps}]{Cyber-Physical Systems}
\acrodef{nsf}[\textsc{nsf}]{National Science Foundation}
\acrodef{pi}[\textsc{pi}]{Principal Investigator}
\acrodef{copi}[\textsc{c}o-\textsc{pi}]{Co-Investigator}
\acrodef{gps}[\textsc{gps}]{Global Positioning System}
\acrodef{ud}[\textsc{ud}]{University of Delaware}
\acrodef{sbu}[\textsc{sbu}]{Stony Brook University}
\acrodef{fame}[\textsc{fame}]{Forum to Advance Minorities in Engineering}
\acrodef{hmi}[\textsc{hmi}]{human-machine interaction}
\acrodef{cav}[\textsc{cav}]{connected and automated vehicle}
\acrodef{sst}[\textsc{sst}]{streaming string transducer}
\acrodef{mpc}[\textsc{mpc}]{model predictive control}
\acrodef{od}[\textsc{od}]{origin-destination}
\acrodef{udssc}[\textsc{udssc}]{University of Delaware scaled smart city}
\acrodef{mso}[\textsc{mso}]{monadic second-order}

\let\marginpar\oldmarginpar
\begin{document}
	\graphicspath{{Figures/}}

	
	\titleformat{\section}
	{
		\large\normalfont\rmfamily\bfseries}
	{\thesection.}{.5em}{}
	
	
	
	
	\titleformat{\subsection}
	{\normalfont\rmfamily\fontseries{b}\selectfont\filright\bfseries}
	{\thesubsection.}{.5em}{}

	

	\thispagestyle{plain}
	
%
	
	
	\thispagestyle{plain}
	\thispagestyle{fancy}
	\renewcommand{\headrulewidth}{0pt}
	\chead{\vspace{0.05in} \bfseries Addressing Mixed Traffic Through Platooning of Vehicles\\ Andreas A. Malikopoulos\\Cornell University}
	\cfoot{\thepage}
	
	\setcounter{page}{1}


\label{section:introduction}
\subsection*{Abstract}
Connected and automated vehicles (CAVs) provide the most intriguing opportunity for enabling users to better monitor transportation network conditions and make better operating decisions to improve safety and reduce pollution, energy consumption, and travel delays. 
While several studies have shown the benefits of CAVs in reducing energy and alleviating traffic congestion in specific traffic scenarios, e.g., crossing signal-free intersections, merging at roadways and roundabouts, cruising in congested traffic, passing through speed reduction zones, and lane-merging or passing maneuvers, most of these efforts have focused on $100$\% CAV penetration rates without considering human-driven vehicles (HDVs). One key question that still remains unanswered is ``how can CAVs and HDVs be coordinated safely?'' 
In this paper, we report on an optimal control framework to coordinate CAVs and HDVs at any  traffic scenario.
The idea is to directly control the CAVs to force the trailing HDVs to form platoons. Thus, we indirectly control the HDVs by coordinating the platoon of HDVs led by CAVs.

\subsection*{Introduction}
\label{STA}

\begin{wrapfigure}{r}{0.7\textwidth}
	\vspace{-2ex}
	\centering
	\includegraphics
	[width=0.7\textwidth]
	{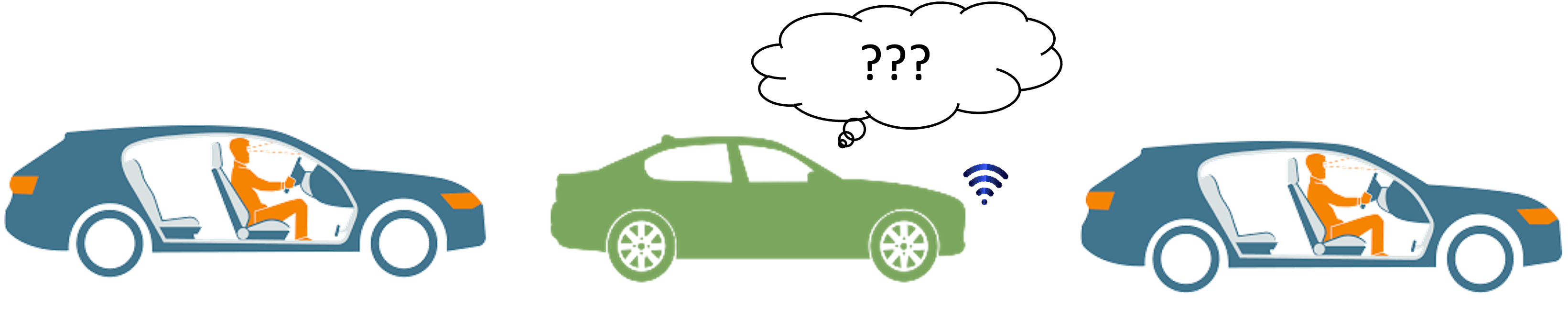}
	\caption{\small
		A mixed traffic environment consisting of human-driven and connected automated vehicles.
	}
	\label{fig:1}
	\vspace{-2ex}
\end{wrapfigure}

The existence of data and shared information in connected and automated vehicles (CAVs) is associated with significant technical challenges and gives rise to a new level of complexity \cite{Malikopoulos2016b} in modeling and control \cite{Ferrara2018} for emerging mobility systems \cite{Malikopoulos2011,zhao2019enhanced} with non-classical informational structures \cite{Malikopoulos2021,Dave2020a,Dave2021a}. Using data and shared information through vehicle-to-vehicle (V2V) and vehicle-to-infrastructure (V2I) communication, we can aim at optimally controlling  CAVs at different traffic scenarios, e.g., crossing signal-free intersections, merging at roadways and roundabouts, cruising in congested traffic, passing through speed reduction zones, and lane-merging or passing maneuvers. These scenarios, along with the driver's responses to unexpected events, contribute to traffic congestion \cite{Malikopoulos2013}.
It is expected that CAVs will gradually penetrate the market, interact with human-driven vehicles (HDVs), and contend with V2V and V2I communication limitations, e.g., bandwidth, dropouts, errors, and delays, as well as perception delays, lack of state information, etc. However, in a mixed traffic environment consisting of CAVs and HDVs (Fig. \ref{fig:1}), HDVs pose significant modeling and control challenges to the CAVs due to the stochastic nature of the human-driving behavior. Therefore, different penetration rates of CAVs can significantly alter transportation efficiency and safety. 
As we move to increasingly diverse mobility systems with different penetration rates of CAVs, new approaches are needed to optimize the impact on system behavior of the interplay between CAVs and HDVs at different traffic scenarios.  While several studies have shown the benefits of CAVs in reducing energy and alleviating traffic congestion in specific traffic scenarios, most of these efforts have focused on $100$\% CAV penetration rates without considering HDVs. 

The objective of this paper is to report on an optimal control framework aimed at coordinating CAVs and HDVs in any  traffic scenario, e.g., crossing signal-free intersections, merging at roadways and roundabouts, cruising in congested traffic, passing through speed reduction zones, and lane-merging or passing maneuvers.
Our framework aims to directly control the CAVs to force the trailing HDVs to form platoons. Thus, we indirectly control the HDVs by coordinating the platoon of HDVs led by CAVs.
In this context, we address the problem of  platoon formation of HDVs  by only controlling the CAVs in the transportation network. 
 
There have been several research efforts reported in the literature addressing the problem \cite{schwarting2018planning}  of the interactions between AVs with HDVs using probabilistic computation tree logic \cite{sadigh2014data}, social value orientation \cite{schwarting2019social,buckman2019sharing,Le2022CDC}, classification on the drivers' behavior \cite{markkula2020defining}, predictions of vehicle trajectories \cite{chandra2020forecasting,wang2021socially,song2022learning}, reinforcement learning (RL) \cite{saxena2020driving}, multi-agent RL inverse RL  \cite{sadigh2016planning,sadigh2018planning}, hierarchical trajectory planning games \cite{fisac2019hierarchical}, multi-agent forecasting  \cite{rhinehart2019precog}.
Other studies have focused on quantifying the impact of AVs on energy and efficiency \cite{Martin2011,Martin2011a,firnkorn2015,Martin2016,Chen2016} while other studies \cite{Martin2011,Martin2011a,Martin2016} have shown that AVs can  decrease  emissions \cite{Martin2011a}. There have also been efforts in investigating the feasibility and potential environmental impacts \cite{Shaheen2014} of shared AVs \cite{chong2013a,Ford2012,Rigole2014,Bischoff2016,Dia2017,Merlin2017,Metz2018,Fiedler2018,Lu2018}. Some authors have focused on the cost-benefit analysis of a mobility system with  AVs \cite{Burns2012,Fagnant2014,Bosch2017,Moorthy2017}, while some other studies have investigated its impact on vehicle ownership by using surveys or comparable analyses with conventional car-sharing systems \cite{Truong2017,tussyadiah2017,Foldes2018,Hao2018,Menon2018}. Several survey papers are providing a good review of related topics (e.g., see \cite{Jorge2013,Agatz2012a,Furuhata2013a,Brandstatter2016,lavieri2017,jittrapirom2017,utriainen2018}).

In a typical commute, we encounter traffic scenarios that include crossing intersections, merging at roadways and roundabouts,  cruising in congested traffic, passing through speed reduction zones, and lane-merging or passing maneuvers.
In an emerging mobility system with a 100\% penetration rate of CAVs, i.e., without HDVs, we can coordinate the CAVs in such scenarios to alleviate congestion, eliminating stop-and-go driving and thus improving energy, GHG emissions, and travel time. Several research efforts have been reported in the literature towards developing control algorithms for coordinating CAVs to alleviate congestion. In 2004, Dresner and Stone \cite{Dresner2004} proposed the use of the reservation scheme to control a single intersection of two roads. In their approach, each vehicle requests the reservation of the space-time cells to cross the intersection at a particular time interval defined from the estimated arrival time to the intersection. 
Since then, numerous reservation and scheduling approaches have been reported in the literature to achieve safe and efficient control of traffic through intersections \cite{DeLaFortelle2010, Dresner2008,Colombo2015,Ahn2014,DeCampos2015a,Ahn2016,fayazi2018mixed,guney2020scheduling}.
Colombo and Del Vecchio \cite{Colombo2015} designed an intersection controller for vehicles that only intervenes and overrides the driver's control action when necessary. These results have been extended to include uncontrolled human drivers \cite{Ahn2014}. In a sequel paper, Ahn and Del Vecchio \cite{Ahn2016} solved the supervisory problem for the first-order dynamics without considering the rear-end collision avoidance constraint using a mixed-integer program.

\begin{wrapfigure}{r}{0.55\textwidth}
	\vspace{-2ex}
	\centering
	\includegraphics
	[width=0.55\textwidth]
	{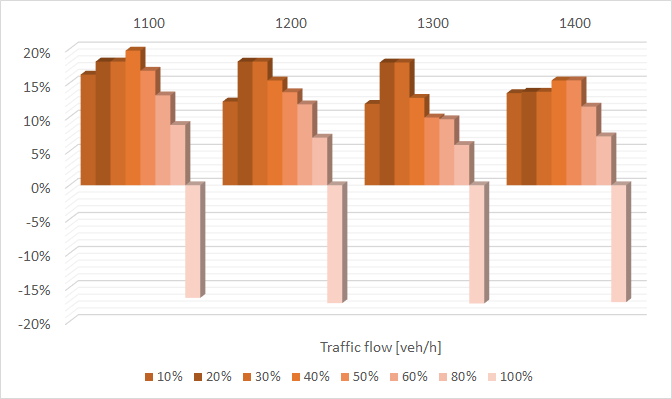}
	\caption{\small
		Fuel consumption implications for different penetration rates of connected and automated vehicles in a merging at highways on-ramp scenario \cite{Malikopoulos2018d}.
	}
	\label{fig:fuel}
	\vspace{-2ex}
\end{wrapfigure}

There have also been efforts using model predictive control \cite{mirheli2019consensus,hult2018optimal,kamal2014vehicle}, fuzzy logic \cite{Onieva2012}, and navigation functions \cite{Makarem2012} to investigate coordination of CAVs at signal-free intersections. 
Queueing theory has also been used to address this problem by modeling the coordination of CAVs as a polling system with two queues and one server that determines the sequence of times assigned to the vehicles on each road \cite{Miculescu2014}.  Some of the methods presented in the literature have focused on multi-objective optimization problems  \cite{Kamal2013a, Kamal2014, Campos2014, Makarem2013, qian2015}. More recently, a study \cite{Ratti2016} indicated that transitioning from intersections with traffic lights to autonomous intersections, where vehicles can coordinate and cross the intersection without the use of traffic lights, has the potential of doubling capacity and reducing delays.

Over the last several years, a decentralized optimal control framework and algorithms were developed for coordinating in real time CAVs at $100$\% penetration rate at different traffic scenarios, e.g., intersections, merging roadways, roundabouts, and corridors \cite{Malikopoulos2017,Malikopoulos2018c, Malikopoulos2016a,Malikopoulos2020,Mahbub2020automatica, Rios-Torres2,chalaki2020TCST,chalaki2020TITS,chalaki2020hysteretic, Beaver2020DemonstrationCity,chalaki2019optimal,chalaki2020experimental,jang2019simulation,Mahbub2019ACC,mahbub2020ACC-2,mahbub2020sae-1,mahbub2020sae-2,Malikopoulos2018b,malikopoulos2019ACC,Zhao2018CTA,Zhao2018ITSC,Malikopoulos2019CDC,Zhao2019CCTA-1,zhao2019CCTA-2,chalaki2019zero,Ray2021DigitalCity,chalaki2021CSM,ChalakiCBF2022}.
The framework provides a closed-form analytical solution that exists under certain conditions \cite{Mahbub2020ACC-1,mahbub2020Automatica-2}, and which, based on Hamiltonian analysis, yields for each CAV the optimal acceleration/deceleration at any time in the sense of minimizing fuel consumption. 
Several reinforcement learning approaches have been also reported to address these problems \cite{chalaki2020hysteretic,chalaki2020ICCA,jang2019simulation,Sumanth2021,chalaki2021RobustGP}. 
A detailed discussion of the research efforts reported in the literature in this area can be found in \cite{Malikopoulos2016a} and \cite{Guanetti2018}. 

\begin{wrapfigure}{r}{0.55\textwidth}
	\vspace{-2ex}
	\centering
	\includegraphics
	[width=0.55\textwidth]
	{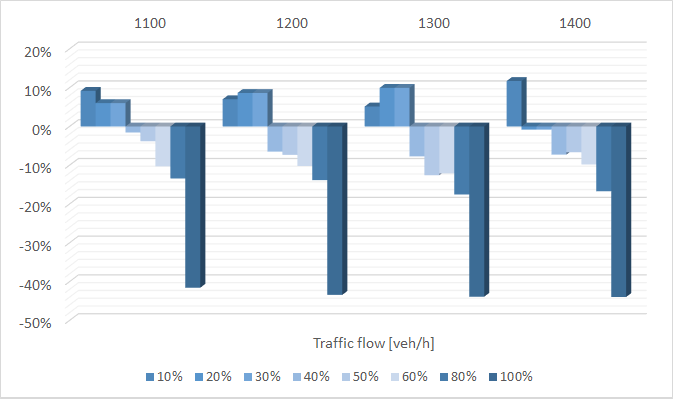}
	\caption{\small
		Travel time implications for different penetration rates of connected and automated vehicles in a merging at highways on-ramp scenario \cite{Malikopoulos2018d}.
	}
	\label{fig:travel}
	\vspace{-2ex}
\end{wrapfigure}

\begin{wrapfigure}{r}{0.55\textwidth}
	\vspace{-2ex}
	\centering
	\includegraphics
	[width=0.55\textwidth]
	{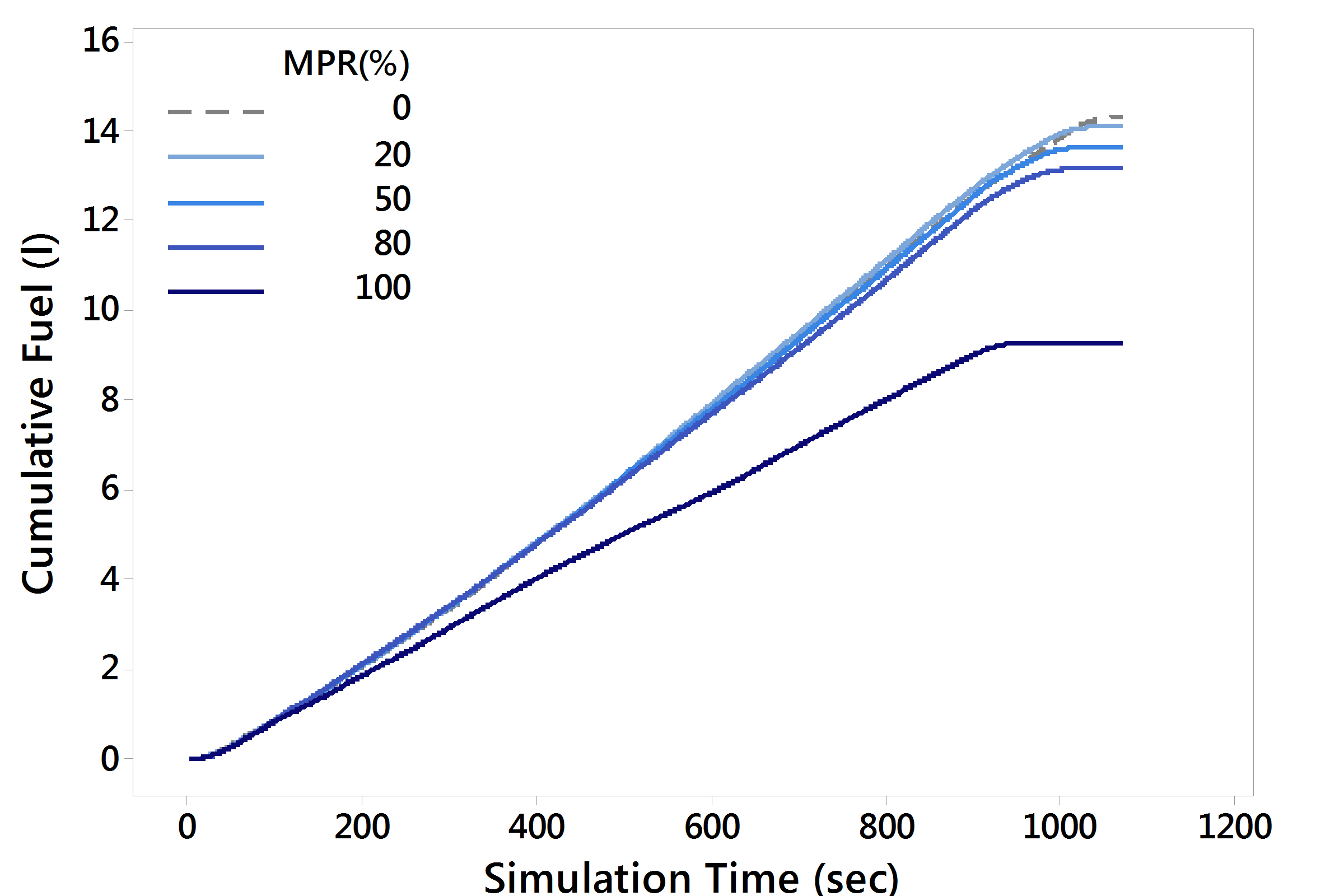}
	\caption{\small
		Impact of different market penetration rates (MPR) of connected and automated vehicles on fuel consumption in merging at roundabouts \cite{Zhao2018CTA}.
	}
	\label{fig:fuel2}
	\vspace{-2ex}
\end{wrapfigure}

In two recent studies \cite{Malikopoulos2018d, Zhao2018CTA}, we investigated the impact of different market penetration rates (MPR) of CAVs, e.g., 0\% to 100\%, in fuel consumption and travel time for two traffic scenarios: (1) merging at highways on-ramp and (2) merging at roundabouts. We observed that as we decrease the penetration rate of CAVs, fuel consumption and travel time deteriorate. Figure \ref{fig:fuel} shows our findings on fuel consumption savings in a merging at highways on-ramp traffic scenario for different traffic volumes and different penetration rates of CAVs. Obviously, only at 100\% penetration rates of CAVs we have any improvements. For travel time (Fig. \ref{fig:travel}), there seem to be improvements up to 40\% penetration rates of CAVs. Figure \ref{fig:fuel2} shows the impact of different MPR of CAVs on fuel consumption in a traffic scenario that includes vehicle merging at roundabouts.  The conclusion that we can draw from these studies is that as we move to increasingly complex  \cite{Malikopoulos2016b} mobility systems  with CAVs and HDVs, we need to develop a system-based approach to improve transportation efficiency in terms of fuel consumption and travel time. 

The  proposed framework advances the state of the art by coordinating CAVs with HDVs at any given traffic scenario, e.g., crossing signal-free intersections, merging
at roadways and roundabouts, cruising in congested traffic, passing through speed reduction zones, and lane-merging or passing maneuvers. 
The novel feature that sharply distinguishes this framework from previous approaches reported in the literature to date is that, in our proposed approach, we directly control the CAVs to force the trailing HDVs to form platoons.
Thus, we indirectly control the HDVs by coordinating the platoon of HDVs led by CAVs.
Recent preliminary results have shown that we can effectively control CAVs to facilitate platoon formation of HDVs  \cite{mahbub2021_platoonMixed,Beaver2021Constraint-DrivenStudy,mahbub2022NHM,mahbub2022_ifac,mahbub2022ACC,mahbub2023_automatica,mahbub2022NHM} and then coordinate the platoons
\cite{Kumaravel:2021uk,Kumaravel:2021wi}  to improve transportation efficiency.

\section{Control Framework} 
\label{section:research}

CAVs are typical cyber-physical systems \cite{Malikopoulos2022a} where the cyber component (e.g., data and shared information through V2V and V2I communication) can aim at optimally controlling the physical entities (e.g., CAVs, non-CAVs) \cite{Cassandras2017}. The cyber-physical nature of such systems is associated with significant control challenges and gives rise to a new level of complexity in modeling and control. 
As we move to increasingly {diverse} emerging mobility systems (e.g., CAVs, shared mobility) with different penetration rates of CAVs, new approaches are needed to optimize the impact on system behavior of the interplay between CAVs and HDVs at different traffic scenarios.


\subsection{Forming Platoon of HDVs led by CAVs}
In our approach, we  directly control the CAVs to force the trailing HDVs to form platoons, and thus, we indirectly control the HDVs. In this context, we address the problem of vehicle platoon formation in a mixed-traffic environment by only controlling the CAVs. 
We seek to derive an analytical solution of the control input of CAVs along with the conditions under which the solution is feasible, thus  establishing a rigorous control framework that enables platoon formation in a mixed traffic environment with associated boundary conditions.

\begin{wrapfigure}{r}{0.6\textwidth}
	\vspace{-2ex}
	\centering
	\includegraphics
	[width=0.6\textwidth]
	{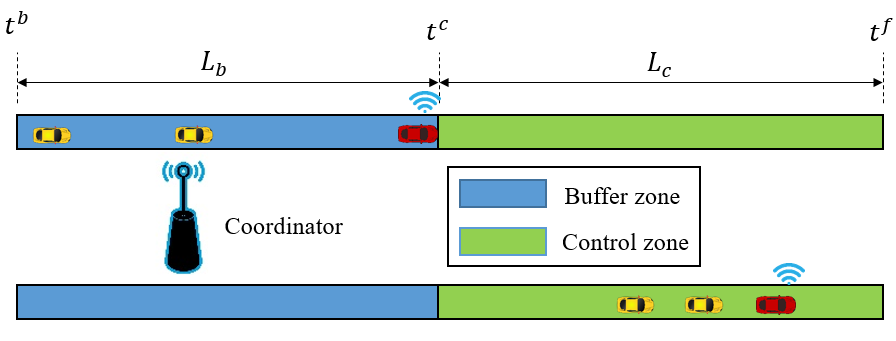}
	\caption{\small
		A CAV (red) traveling with two trailing HDVs (yellow).
	}
	\label{fig:platoon_zone}
	\vspace{-2ex}
\end{wrapfigure}

We consider a CAV followed by one or multiple HDVs traveling in a single-lane roadway of length $L\in \mathbb{R}^+$. We subdivide the roadway into a {buffer zone} of length $L_{b}\in \mathbb{R}^+$, inside of which the HDVs' state information is estimated through roadside units (Fig. \ref{fig:platoon_zone}) (top), and a {control zone} of length $L_{c}\in \mathbb{R}^+$ such that $L=L_{b}+L_{c}$, where the CAV is controlled to form a platoon with the trailing HDVs, as shown in Fig. \ref{fig:platoon_zone} (bottom). The time that a CAV enters the buffer zone, the control zone, and exits the control zone is $t^b, t^c, t^f\in \mathbb{R}^+$, respectively.

Let $\mathcal{N}=\{1,\ldots, N\}$, where $N\in \mathbb{N}$ is the total number of vehicles traveling within the buffer zone at time $t=t^c$, be the set of vehicles considered to form a platoon. 
Here, the leading vehicle indexed by $1$ is the CAV, and the other trailing vehicles in $\mathcal{N}\setminus\{1\}$ are HDVs.  
We denote the set of the HDVs following the CAV to be $\mathcal{N}_{\text{HDV}}=\{2,\ldots, N\}$.
Since the HDVs do not share their local state information with any other vehicle, we consider the presence of a {coordinator} that gathers the state information of the trailing HDVs traveling within the buffer zone. The coordinator, which can be a group of loop-detectors or comparable sensory devices, then transmits the HDV state information to the CAV at each time instance $t\in[t^b, t^c]$ using a standard vehicle-to-infrastructure communication protocol.  At the entry of the control zone at time $t^c$, the CAV uses the state information from the HDVs to derive and apply the control input so that a platoon is formed with the following HDVs at some time $t> t^c$ within the control zone.
The objective of the CAV $1$ is to derive and implement a control input (acceleration/deceleration) at time $t^c \in \mathbb{R}^+$ so that the platoon formation with trailing HDVs in $\mathcal{N}_{\text{HDV}}$ is completed within the control zone at a given time $t^p\in (t^c, t^f]$.

In our framework, we model the longitudinal dynamics of each vehicle $i\in\mathcal{N}$ as a  double-integrator,
\begin{gather}\label{eq:dynamics_pv}
	\dot{p}_i(t) = v_i(t), \quad 
	\dot{v}_i(t) = u_i(t), \quad t\in\mathbb{R}^+,
\end{gather}
where $p_i(t)\in \mathcal{P}_i$, $v_i(t)\in \mathcal{V}_i$, and $u_i(t)\in\mathcal{U}_i$ are the position of the front bumper, speed, and control input (acceleration/deceleration) of vehicle $i\in \mathcal{N}$.
The reason that we model each vehicle as a double integrator is to derive an optimal, closed-form, analytical solution of the control input (acceleration/deceleration) of the CAV in real time, which would be, otherwise, rather infeasible to do so. However, to implement the optimal control input in the vehicle, a PI or PID controller designed to capture the vehicle dynamics is used to track the optimal acceleration/deceleration given by the model in \eqref{eq:dynamics_pv}. So, the optimal control problem is eventually ``translated" into a tracking control problem (tracking the optimal acceleration/deceleration).

Let $\textbf{x}_{i}(t)=\left[p_{i}(t) ~ v_{i}(t)\right] ^{T}$ denote the state vector of each vehicle $i\in\mathcal{N}$, taking values in $\mathcal{X}_{i}%
=\mathcal{P}_{i}\times\mathcal{V}_{i}$. 
The state space $\mathcal{X}_{i}$ for each vehicle $i\in\mathcal{N}$ is closed concerning the induced topology on $\mathcal{P}_{i}\times \mathcal{V}_{i}$, and thus, it is compact.
The speed $v_i(t)$ and control input $u_i(t)$ of each vehicle $i\in \mathcal{N}$ are subject to the following constraints,
\begin{align}\label{eq:state_control_constraints}
	0< v_{\min} \le v_i(t)\le v_{\max}, \quad u_{i,\min} \le u_i(t) \le u_{i,\max}, \quad t\in\mathbb{R}^+,
\end{align}
where $v_{\min}$ and $v_{\max}$ are the minimum and maximum allowable speed of the considered roadway, respectively, and $u_{i,\min}$ and $u_{i,\max}$ are the minimum and maximum control input of each vehicle $i\in \mathcal{N}$, respectively.

\begin{wrapfigure}{r}{0.5\textwidth}
	\vspace{-2ex}
	\centering
	\includegraphics
	[width=0.5\textwidth]
	{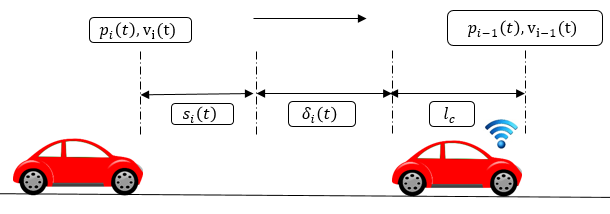}
	\caption{\small
		Predecessor-follower coupled car-following dynamic.
	}
	\label{fig:problem_formulation}
	\vspace{-2ex}
\end{wrapfigure}

The dynamics \eqref{eq:dynamics_pv} of each vehicle $i\in\mathcal{N}$ can take different forms based on the consideration of connectivity and automation.
For the CAV $1\in\mathcal{N}$, the control input $u_1(t)$ can be derived and implemented within the control zone. We introduce and discuss the structure of the control zone in detail in Task 1.2.
The dynamic following spacing $s_i(t)$ between two consecutive vehicles $i \text{ and }(i-1)\in\mathcal{N}$ is, ${ s_i(t)= \rho_i\cdot v_i(t)+ s_0,}$
where  $\rho_i$ denotes a desired time gap that each HDV $i\in\mathcal{N}_{\text{HDV}}$ maintains while following the preceding vehicle, and $s_0$ is the standstill distance denoting the minimum bumper-to-bumper gap at stop.
The {platoon gap} $\delta_i(t)$ is the difference between the bumper-to-bumper inter-vehicle spacing and the dynamic following spacing {$s_i(t)$} (see Fig. \ref{fig:problem_formulation}) between two consecutive vehicles $i \text{ and }(i-1)\in\mathcal{N}$, i.e., $\delta_i(t)=p_{i-1}(t)- p_i(t)-s_i(t)-l_c,$
where $l_c$ is the length of each vehicle $i\in\mathcal{N}$.

In our approach, we adopt the optimal velocity car-following model \cite{bando1995dynamical} to define the predecessor-follower coupled dynamics (see Fig. \ref{fig:problem_formulation}) of each HDV $i\in \mathcal{N}_{\text{HDV}}$ as follows,
\begin{gather}\label{eq:hdv_dynamics}
	{{u}_i(t) = \alpha (V_i(\delta_i(t-\eta_i),s_i(t-\eta_i)) -v_i(t-\eta_i)),}
\end{gather}
where $\alpha$ denotes the control gain representing the driver's sensitivity coefficient, {$\eta_i$ is the driver's perception delay with a known upper bound $\bar{\eta}$}, and $V_i(\delta_i(t),s_i(t))$ denotes the equilibrium speed-spacing function,
\begin{gather}
	V_i(\delta_i(t),s_i(t))=
	\begin{array}
		[c]{ll}%
		{\frac{v_{\max}}{2}(\tanh(\delta_i(t))}{+\tanh(s_i(t))).}
	\end{array}
	\label{eq:V(s)}
\end{gather}

Based on \eqref{eq:V(s)}, the driving behavior of each HDV $i\in\mathcal{N}_{\text{HDV}}$ depends on two different modes; (a) {decoupled free-flow mode}: when $\delta_i(t)>0$, each HDV converges to the maximum allowable speed $v_{\max}$, and cruises through the roadway decoupled from the state of the preceding vehicle, 
and (b) {coupled following mode}: when $ \delta_i(t)\le 0$, the HDV dynamics becomes coupled with the state of the preceding vehicle $(i-1)\in\mathcal{N}$, and $v_i(t)$ converges to $v_{i-1}(t)$. Note that if there is no preceding vehicle, we set $\delta_i(t)=\infty$ that activates the decoupled free-flow mode, which results in $v_i(t)$ converging to $v_{\max}$.

The car-following model \eqref{eq:hdv_dynamics} is platoon-stable, i.e., bounded speed fluctuation between two consecutive vehicles in coupled following mode decays exponentially as time progresses \cite{bando1995dynamical}.
Next, we define the information set $\mathcal{I}_1(t)$ of the CAV $1\in\mathcal{N}$ which has the following structure, $\mathcal{I}_1(t) = \{\textbf{x}_1(t), \textbf{x}_{2:N}(t)\},
		~ t\in[t^b,t^c],$	where $\textbf{x}_{2:N}(t)=[\textbf{x}_2(t),\ldots, \textbf{x}_N(t)]^T$.
The steady-state traffic flow between two consecutive vehicles $i\text{ and }(i-1)\in\mathcal{N}$ are established if the platoon gap $\delta_i(t)$ does not vary with time, and {speed fluctuation $\Delta v_i(t):=v_i(t)-v_{i-1}(t)$ is zero} \cite{rothery1992car}, i.e., 
$\delta_i(t)= c_i,~c_i\in \mathbb{R}, \text{ and }{\Delta v_i(t) =0}.$

Next, we formulate the problem of platoon formation in a mixed environment as follows:
\begin{problem}\label{prob:1}
	Given the information set $\mathcal{I}_1(t)$ at time $t=t^c$, the objective of the CAV $1\in\mathcal{N}$ is to derive the control input $u_1(t)$ so that the HDVs in ${N}_{HDV}$ are forced to form a platoon at some time $t^p\in (t^c, t^f]$ within the control zone while the following conditions hold,
	\begin{align}
		v_i(t) = v_{eq},~  &\delta_i(t)=c_i,~c_i\le 0, \quad \forall t\ge t^p,~ \forall i\in\mathcal{N}, \nonumber\\
		&\text{ subject to: }{ \eqref{eq:state_control_constraints},~ p_1(t^p)\le L_c,} \label{eq:platoon_cond_1}
	\end{align}
	where, $v_{eq}$ denotes the equilibrium platoon speed.
\end{problem}

In our problem formulation, we impose the restriction that at $t=t^c$, there exists at least one HDV $i\in\mathcal{N}_{\text{HDV}}$ such that $\delta_i(t^c)>0$. To simplify the formulation and without loss of generality, we consider that $\delta_N(t^c)>0$. This ensures that we do not have the trivial case where the group of vehicles in $\mathcal{N}$ has already formed a platoon at $t=t^c$. 
In the modeling framework presented above, we impose the following {assumption:}
The CAV is on a decoupled free-flow mode while all vehicles have reached steady-state traffic flow  within $[t^b, t^c]$.
Note that we restrict the control of the CAV $1$ only within the control zone so that we have a finite control horizon $[t^c,t^f]$. Outside the control zone, the CAV dynamics follow the car-following model in \eqref{eq:hdv_dynamics}.

We address Problem \ref{prob:1} considering only two vehicles first, i.e., $N=2$ and then generalize the analysis for multiple vehicles, i.e., $N>2$. When the CAV $1\in\mathcal{N}$ applies a control input $u_1(t),~t\in[t^c, t^p]$ to form a platoon with the HDV $2\in\mathcal{N}_{\text{HDV}}$ at time $t=t^p$, two sequential steps take place, namely, (i) the {platoon transition} step, where the HDV $2$ transitions from the decoupled free-flow mode to the coupled following mode at time $t=t^s, ~t^c< t^s< t^p$ such that $\delta_2(t^s)=0$, and (ii) the {platoon stabilization} step, where $v_2(t)$ converges to $v_1(t)$ at time $t=t^p$ such that \eqref{eq:platoon_cond_1} is satisfied, and the platoon becomes stable.

The platoon transition duration $\tau^t$ is the time required for the completion of the platoon transition step, i.e., $\tau^t=t^s-t^c$, and the platoon stabilization duration $\tau^s$ is the time required for the completion of the platoon stabilization step, i.e., $\tau^s=t^p-t^s$. Hence, we have $t^p=t^c+\tau^{t}+\tau^{s}$.
The platoon stabilization duration is $\tau^s=\eta_i+\tau^r$, where $\eta_i$ is the perception delay of HDV $i\in\mathcal{N}_{\text{HDV}}$, and $\tau^r$ is the response time of $\eqref{eq:hdv_dynamics}$ which depends on the driver’s sensitivity coefficient $\alpha$, maximum allowable speed fluctuation, and the choice of equilibrium speed-spacing function in \eqref{eq:V(s)}, and can be computed using standard stability analysis \cite{bando1995dynamical,wilson2011car,zhang2021improved}. Note that, for $N\ge 2$,  additional nonlinearities may impact the computation of $\tau^s$.
In our formulation, we incorporate the upper bound of the perception delay $\bar{\eta}$ to achieve {robustness} such that $\tau^s=\bar{\eta}+\tau^r$, and consider that $\tau^r$ is given a priori. Thus, we focus only on analyzing the platoon transition time $\tau^t$.
We  construct the structure of the control input $u_1(t)$ for the CAV $1\in\mathcal{N}$ for generating a platoon with the trailing HDV $2\in\mathcal{N}$ at time $t^p\in(t^c, t^f]$,
\begin{gather}\label{eq:control_structure}
	u_{1}(t) = \left\{
	\begin{array}
		[c]{ll}%
		u_p, ~u_p\in[u_{\min},0), & t\in[t^c, t^s],\\
		0, ~  & t\in( t^s,t^f].
	\end{array}
	\right.
\end{gather}
Based on \eqref{eq:control_structure}, the realization of the control input $u_1(t)$ of the CAV $1\in\mathcal{N}$, where $u_1(t)=u_p\in(0,u_{\min}]$ in $t\in[t^c, t^s]$, yields a linearly decreasing $v_1(t)$ in $t\in[t^c, t^s]$.
We  investigate whether there will always exist an unconstrained control input parameter $u_p$ in \eqref{eq:control_structure} such that a vehicle platoon can be formed with HDV $2\in\mathcal{N}$ at time $t=t^p$ given the following condition,
	\begin{equation}\label{eq:u_t_relation_1}
		2\delta_2(t^c) + u_p\cdot(\tau^t)^2 = 0.
	\end{equation}
The advantage of condition \eqref{eq:u_t_relation_1} is that it provides an intuitive result. Namely, as $u_p\rightarrow 0$, we have $\tau^t \rightarrow \infty$, which implies that platoon formation will never occur. If $u_p>0$, then \eqref{eq:u_t_relation_1} yields an infeasible $\tau^t$. Therefore, $u_p$ has to be strictly negative for platoon formation, i.e., CAV will need to slow down. Note, from \eqref{eq:u_t_relation_1}, for $\delta_2(t)>0$ and $t\in \mathbb{R}^+$, we have $u_p<0$.

\subsection{Feasibility of the Platoon Formation Time}

Our analysis does not explicitly incorporate the state and control constraints in \eqref{eq:state_control_constraints} and the terminal constraint in \eqref{eq:platoon_cond_1}. For a given platoon formation time $t^p$, the corresponding control input derived from \eqref{eq:u_t_relation_1} can violate constraints in \eqref{eq:state_control_constraints}. We explore the existence of a feasible region of $\tau^t$ that yields an admissible control input parameter $u_p$ in \eqref{eq:u_t_relation_1}.

Suppose that for CAV $1\in\mathcal{N}$, a platoon transition duration $\tau^{t}$ has associated control input parameter $u_p$ derived from \eqref{eq:u_t_relation_1}.
Using \eqref{eq:dynamics_pv}, we have,
$v_1(t^c+\tau^{t})=v_1(t^c)+u_p \cdot \tau^{t}.$ Since $u_p\in[u_{\min},0)$, we require that $v_1(t^c+\tau^t)\ge v_{\min}$ to satisfy the state constraint in \eqref{eq:state_control_constraints}. Substituting $v_1(t^c+\tau^t)$ in the above inequality, we get $ u_p\cdot \tau^t \ge v_{\min}-v_1(t^c)$. Finally, substituting $u_p$ from \eqref{eq:u_t_relation_1} in the last equation we have
	\begin{gather}
	\tau^t \ge \max\bigg\{ \bigg(\frac{-2\delta_2(t^c)}{u_{\min}}\bigg)^{\frac{1}{2}}, \frac{2\delta_2(t^c)}{v_1(t^c)-v_{\min}} \bigg\}. \label{eq:lem:3}
\end{gather}
The minimum speed value $v_{\min}$ in \eqref{eq:lem:3} indicates the allowable speed perturbation during the platoon stabilization step. Hence, $v_{\min}$ should be selected appropriately to ensure the local stability of the platoon \cite{wilson2011car, zhang2021improved}.
We need to investigate the conditions under which, for a given platoon formation time $t^p$ and platoon stabilization duration $\tau^s$, the corresponding platoon transition duration $\tau^t$ is feasible.

For $N>2$, the CAV $1\in\mathcal{N}$ trailed by multiple HDVs $j\in\mathcal{N}_{\text{HDV}}$ and given $\mathcal{I}_1(t^c)$, we have the following conditions, $v_1(t^c)=v_j(t^c)=v_{\max}$ for all $j\in\mathcal{N}_{\text{HDV}}$, and there exists $j\in\mathcal{N}_{\text{HDV}}$ such that $\delta_j(t^c)>0$.
For a CAV $1\in\mathcal{N}$ followed by $N \in \mathcal{N_{\text{HDV}}}$ HDVs, the cumulative platoon gap $\Delta(t)$ at time $t\in[t^c,t^f]$ is, $\Delta(t) = p_1(t)-p_{N}(t)-\sum_{j=2}^{N} (s_j(t)+ l_c).$

\subsection{Coordinating Platoon of Vehicles led by CAVs}

After the platoons of HDVs led by CAVs are formed within the {buffer zone} (Fig. \ref{fig:platoon_zone}), the next step is to coordinate these platoons in a traffic scenario.
Next, we provide a framework for coordinating platoons of HDVs led by CAVs at different traffic scenarios, e.g.,  merging at roadways and roundabouts, cruising in congested traffic, passing through speed reduction zones, and lane-merging or passing maneuvers. Although this framework can be applied to any  traffic scenario, we use an intersection as a reference to provide the fundamental ideas. 
The proposed framework includes {(1)} an upper-level optimization that yields for each platoon its optimal time trajectory and lane to pass through a given traffic scenario while alleviating congestion; and {(2)} a low-level optimization that yields for each leader of the platoon, i.e., the CAV, its optimal control input (acceleration/deceleration). 

\begin{wrapfigure}{r}{0.4\textwidth}
	\vspace{-2ex}
	\centering
	\includegraphics
	[width=0.4\textwidth]
	{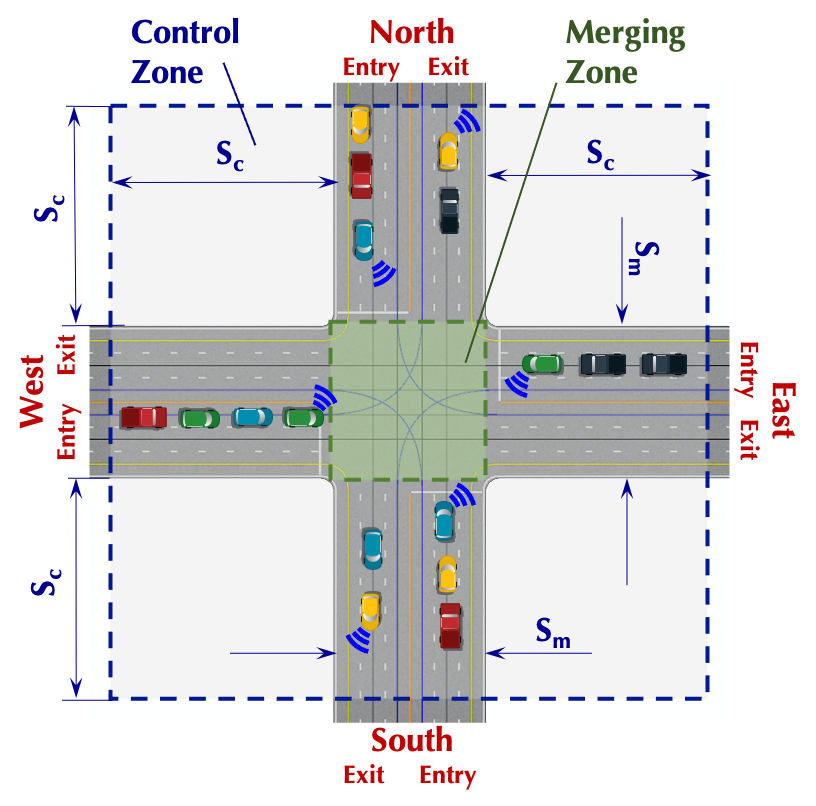}
	\caption{\small
		A signal-free intersection with human-driven led by connected automated vehicles.
	}
	\label{fig:6}
	\vspace{-2ex}
\end{wrapfigure}

We consider platoons of HDVs led by CAVs crossing a signal-free intersection (Fig. \ref{fig:6}). The region at the center of the intersection, called {merging zone}, is the area of potential lateral collision of  the vehicles. The intersection has a {control zone} inside of which the CAVs, which lead the platoons, can communicate with each other. The distance from the entry of  the control zone until the entry of the merging zone is $S_c$ and although it is not restrictive, we consider it to be the same for all entry points of the control zone. We also consider the merging zone to be a square of side $S_m$ (Fig. \ref{fig:6}). Note that  $S_c$ could be in the order of hundreds of meters depending on the CAVs' communication range capability, while $S_m$ is the length of a typical intersection. The aforementioned values of the intersection's geometry are not restrictive in our modeling framework and are used only to determine the total distance traveled by each vehicle inside the control zone.

Let $\mathcal{N}(t)=\{1,\ldots,N(t)\}$, $N(t)\in\mathbb{N}$, be the set of CAVs inside the control zone at time $t\in\mathbb{R}^{+}$. Let $t_{i}^{f}$ be the time for CAV $i$ to exit the control zone.
There are a number of ways to determine $t_{i}^{f}$ for each CAV $i$. The policy, which  determines the time $t_{i}^{f}$ that each CAV $i$ exits the control zone,  is the
result of the upper-level optimization problem, which can aim at maximizing the throughput at the intersection. On the other hand, deriving the optimal control input (minimum acceleration/deceleration) for each CAV $i\in\mathcal{N}(t)$ from the time $t_{i}^{0}$ it enters the control zone until the target $t_{i}^{f}$ is the result of a low-level optimization problem that can aim at minimizing the energy of each CAV, and thus of their entire platoon.

In what follows, we propose a two-level, joint optimization framework: (1) an upper-level optimization that yields for each CAV $i\in\mathcal{N}(t)$, with a given origin (entry of the control zone) and desired destination (exit of the control zone), (a) the minimum time $t_{i}^{f}$ to exit the control zone and (b) optimal path including the lanes that each platoon should be occupying while traveling inside the control zone; and (2) a low-level optimization  that yields, for CAV $i\in\mathcal{N}(t)$, its optimal control input (acceleration/deceleration) to achieve the optimal path and $t_{i}^{f}$ derived in (1) subject to the state, control, and safety constraints.
The two-level optimization framework is executed by each CAV $i\in\mathcal{N}(t)$ as follows. When a CAV $i$ enters the control zone at $t_{i}^{0}$, it communicates with the coordinator  and receives information about the path of each platoon cruising inside the control zone. 
Then, CAV $i$ solves the upper-level optimization problem and derives the minimum time $t_{i}^{f}$ to exit the control zone along with the appropriate lanes that should occupy. Once  CAV $i$ derives $t_{i}^{f}$, then it sends it to the coordinator, who  estimates the time $t_{j}^{f}$, $j\in\mathcal{N}_{\text{HDV}},$ that the last vehicle in the platoon led by CAV $i$ exits the merging zone. The coordinator shares this information with the next CAV entering the control zone, repeating the process.
The minimum time $t_{i}^{f}$ of the upper-level optimization problem is the input of the low-level optimization problem. 
The implications of the proposed optimization framework are that the platoons do not have to come to a complete stop at the intersection, thereby conserving momentum and fuel while also improving travel time. Moreover, by minimizing the control input (acceleration/deceleration) of CAV $i$, essentially, we minimize the acceleration/deceleration of the entire platoon, and thus, we have direct benefits in fuel consumption \cite{Malikopoulos2008b}.

Before we formulate the low-level optimization problem, we consider that the solution of the upper-level optimization problem is given, and thus, the minimum time $t_{i}^{f}$ for each CAV $i\in\mathcal{N}(t)$ is known. 
We first focus on a low-level optimization problem that yields for each CAV $i$ the optimal control input (acceleration/deceleration) to achieve the minimum time $t_{i}^{f}$ subject to the state, control, and safety constraints. The low-level optimization problem is formulated as follows:

\begin{problem} \label{problem1}
	If $t_{i}^{f}$ is determined, the low-level optimal control problem for each CAV $i\in\mathcal{N}(t)$ is to minimize the cost functional $J_{i}(u(t))$, which is the $L^2$-norm of the control input in $[t_i^0, t_i^f]$, i.e.,
	\begin{gather}\label{eq:decentral}
		\min_{u_i(t)\in \mathcal{U}_i} J_{i}(u_i(t)) = \frac{1}{2} \int_{t^0_i}^{t^f_i} u^2_i(t)~dt,
	\end{gather}
\end{problem}
given $t_{i}^{0}\text{, }v_{i}^{0}\text{, } p_{i}(t_{i}^{0})\text{, }t_{i}^{f}\text{,
}p_{i}(t_{i}^{f}),$ where $p_{i}(t_{i}^{0})=0$ and $p_{i}(t_{i}^{f})=p_i^f$. 
Problem \eqref{problem1} is subject to the following {constraints: (1)} state and control constraints, {(2)} rear-end safety constraint, and {(3)} lateral collision avoidance constraint inside the merging zone.

To derive the analytical solution of Problem \eqref{problem1}, we  use Hamiltonian analysis. We first start with the unconstrained arc. If the solution violates any state or control constraints, then the unconstrained arc is pieced together with the arc corresponding to the violated constraint. The two arcs yield a set of algebraic equations that are solved simultaneously using the boundary conditions and interior constraints between the arcs. If the resulting solution, which includes the determination of the optimal switching time from one arc to the next one, violates another constraint, then the last two arcs are pieced together with the arc corresponding to the newly violated constraint, and we resolve the problem with the three arcs pieced together. The three arcs  yield a new set of algebraic equations that must be solved simultaneously using the boundary conditions and interior constraints between the arcs. The resulting solution includes the optimal switching time from one arc to the next one. The process is repeated until the solution does not violate any other constraints.

\subsection{Upper-level Optimization Problem} 
When a CAV $i\in\mathcal{N}(t)$, enters the control zone, it communicates with the coordinator and solves an upper-level optimization problem. The solution to this problem yields for CAV $i$ the  time trajectory $t_{p_i}(p_i)$ for the CAV and, as a result, for the entire platoon. 
In our exposition, we seek to derive the minimum time $t_{i}^{f}$ that CAV $i$ exits the control zone without activating any of the state and control constraints of the low-level optimization Problem \ref{problem1}. Therefore, the upper-level optimization problem  should yield a $t_{i}^{f}$ such that the solution of the low-level optimization problem will result in the unconstrained case. 
There is an apparent trade-off between the two problems. The lower the value of $t_{i}^{f}$ in the upper-level problem, the higher the value of the control input in $[t_{i}^{0}, t_{i}^{f}]$ in the low-level problem. 
The low-level problem is directly related to minimizing energy for each vehicle (individually optimal solution) in the platoon. On the other hand, the upper-level problem is related to maximizing the throughput of the intersection, thus eliminating stop-and-go driving and travel time (system optimal solution). Therefore, by seeking a solution for the upper-level problem, which guarantees  that none of the state and control constraints becomes active, may be considered an appropriate compromise between the two. 

For each CAV $i\in\mathcal{N}(t)$ the optimal position of the unconstrained arc in the following form \cite{Malikopoulos2020}
\begin{gather}
	p^{*}_{i}(t) = \phi_{i,3} \cdot t^3 +\phi_{i,2} \cdot t^2 + \phi_{i,1} \cdot t +\phi_{i,0} , ~ t\in [t_{i}^{0}, t_{i}^{f}], \label{eq:upper_p}%
\end{gather}
where $\phi_{i,3}\neq 0, \phi_{i,2}, \phi_{i,1}, \phi_{i,0}\in\mathbb{R}$ are the constants of integration that will be derived in the Hamiltonian analysis. 
For each $i\in\mathcal{N}(t)$, the optimal position \eqref{eq:upper_p} is a real-valued continuous and differentiable function $\mathbb{R}_{\ge 0}\mapsto \mathbb{R}_{\ge 0}$. It is also a strictly increasing function with respect to $t\in\mathbb{R}_{\ge 0}$.
The optimal speed and control are given by
\begin{alignat}{3}
	v^{*}_{i}(t) & = 3\phi_{i,3} \cdot t^2 +2\phi_{i,2} \cdot t + \phi_{i,1}, \quad && t \in [t_{i}^{0}, t_{i}^{f}], && \label{eq:upper_v} \\
	u^{*}_{i}(t) & = 6\phi_{i,3} \cdot t + 2\phi_{i,2}, && t \in [t_{i}^{0}, t_{i}^{f}]. && \label{eq:upper_u}
\end{alignat}
For each CAV $i\in\mathcal{N}(t)$, the optimal position $p_i^*(t)$ given by \eqref{eq:upper_p} is a one-to-one function for all $t\in[t_{i}^{0}, t_{i}^{f}]$.
Since for each $i\in\mathcal{N}(t)$, $p_i^*(t)$ is a strictly increasing function with respect to $t\in\mathbb{R}_{\ge 0}$, it follows from the mean value theorem that for all $t_1, t_2 \in [t_{i}^{0}, t_{i}^{f}]$ with $t_1 \neq t_2$, we have $p_i^*(t_{1}) \neq p_i^*(t_{2})$.
Therefore,\eqref{eq:upper_p} is a bijective function and its inverse exists.

We rewrite the cubic polynomial function \eqref{eq:upper_p} as
\begin{equation}\label{eqn:upper_p_modified}
	t^3 + \frac{\phi_{i,2}}{\phi_{i,3}} t^2 + \frac{\phi_{i,1}}{\phi_{i,3}} t + \left( \frac{\phi_{i,0}}{\phi_{i,3}} - \frac{p_i}{\phi_{i,3}} \right) = 0 , ~ t\in [t_{i}^{0}, t_{i}^{f}],
\end{equation}
which then can be reduced by the substitution $t = \tau - \frac{\phi_{i,2}}{3\phi_{i,3}}$ to the normal form
\begin{equation}\label{eqn:depressed_cubic}
	\tau^3 + \omega_{i,0} \tau + \left( \omega_{i,1} + \omega_{i,2} p_i \right) = 0,
\end{equation}
where
\begin{gather}\label{eqn:omega1}
	\omega_{i,0} = \frac{\phi_{i,1}}{\phi_{i,3}} - \frac{1}{3}\left(\frac{\phi_{i,2}}{\phi_{i,3}}\right)^2, 
	\omega_{i,1} = \frac{1}{27}\left[2\left(\frac{\phi_{i,2}}{\phi_{i,3}}\right)^3 - \frac{9 \phi_{i,2} \cdot \phi_{i,1}}{(\phi_{i,3}) ^ 2} \right] + \frac{\phi_{i,0}}{\phi_{i,3}}, 
	\omega_{i,2} = - \frac{1}{\phi_{i,3}}.
\end{gather}
We are interested in deriving the expression for the inverse function of \eqref{eq:upper_p} which can be accomplished by finding the root of \eqref{eqn:depressed_cubic}.
Since, for each $i \in \mathcal{N}(t)$, \eqref{eq:upper_p} is a bijective function, there exists an inverse function $p_i^*(t)^{-1}$. The inverse function, $p_i^*(t)^{-1}$ is the time trajectory $t_{p_i}(p_i)$ that yields the time that CAV $i$ is at the position $p_i$ inside the control zone, which can be derived using the Cardano method for cubic polynomials.

For each CAV $i\in\mathcal{N}(t)$,  we will formulate a constrained optimization problem to yield its optimal path in $[t_i^0, t_i^f]$, and thus, the optimal path for each platoon. We start with the introduction of the {cost function} and proceed with the {equality} and {inequality constraints.}
We seek to derive the minimum time $t_{i}^{f^*}$ that a CAV $i\in\mathcal{N}(t)$ exits the control zone without activating any of the state and control constraints of the low-level optimization Problem \ref{problem1}.
For each CAV $i$, the minimum time $t_{i}^{f^*}$ can be derived by minimizing the time trajectory $t_{p_i}(p_i)$ evaluated at $p_i^f$.  
For any fixed $p_i\in[p_i^0, p_i^f]$ of $i\in\mathcal{N}(t)$, since the time trajectory $t_{p_i}(p_i)$ is a function of $\phi_i$, if we vary the constants  $\phi_i$, the time that  $i$ is at the position $p_i$ changes. Thus, in our analysis, we construct the function $f_i:\Phi_i\to[t_i^0, t_i^f]$, $\Phi_i\subset\mathbb{R}^4$, which evaluates the time trajectory at $p_i^f$ and yields the time that each CAV $i$ is located at $p_i^f$ with respect to the variables $\phi_i$, i.e.,	$f_i(\phi_i) = t_{p_i}(p_i^{f}).$
Therefore, to derive the minimum time $t_{i}^{f^*}$ for a CAV $i$, we seek to minimize $f_i(\phi_i)$, with respect to $\phi_i= (\phi_{i,3}, \phi_{i,2}, \phi_{i,1}, \phi_{i,0})$.

The initial and final conditions at the entry and exit of the control zone, respectively, along with the interior constraint $p_(t_i ^ m) = p_i^m$, at the time $t_i^m$ that CAV $i$ enters the merging zone, designate five equality constraints, i.e., $h_i^{(r)}(\phi_i), r=1,\dots, 5$. By contrast, the state and control constraints, along with the rear-end and lateral collision avoidance constraints, designate seven inequality constraints, i.e.,
$g_i^{(m)}(\phi_i), m=1,\dots, 7$.
Therefore, for each CAV $i\in\mathcal{N}(t)$, the upper-level optimization problem is formulated as follows:
\begin{problem} \label{problem2}
	\begin{align}\label{eq:primal}
		\min_{\phi_i} &~ f_i (\phi_i) \nonumber\\
		\text{subject to:}\quad  \phi_i\in\Phi_i, \quad & h_i^{(r)}(\phi_i)=0,~ r=1,\dots, 5, \quad g_i^{(m)}(\phi_i)\le 0, ~m=1,\dots, 7.
	\end{align}
\end{problem}
It can be shown that there is no duality gap in  Problem \eqref{problem2} \cite{Malikopoulos2020}, hence the solution can be derived in real time using standard numerical algorithms.
The proposed framework is implemented as follows. Every time a CAV $i\in\mathcal{N}(t)$ enters the control zone, it formulates and solves Problem \ref{problem2} by communicating with the coordinator.	The solution yields the optimal time trajectory $t_{p_i}(p_i)$ of CAV $i$, and as a result, the minimum time $t_i^f$ to exit the control zone. Then, the optimal trajectory of the entire platoon led by CAV $i$ is estimated by the coordinator.


Problem \ref{problem2} is solved sequentially by each CAV that enters the control zone. Once the time trajectory of a CAV inside the control zone is derived, then it does not change. By inversing $t_{p_i}(p_i)$ and taking the second derivative, CAV $i$ obtains the optimal control input that corresponds to the unconstrained arc. Therefore, the solution of Problem \ref{problem2}, if it exists, guarantees that none of the state and control constraints becomes active in the low-level optimization. If, however, the solution of Problem \ref{problem2} does not exist, then CAV $i$ selects a feasible $t_i^f$ from the coordinator and follows the analysis of the low-level optimization, which includes piecing together the constrained and unconstrained arcs to derive the optimal control input from $t_i^0$ to $t_i^f$.

\section{Concluding Remarks}
\label{section:evaluation}
In this paper, we discussed a control framework aimed at forming platoon under a mixed traffic environment. In this framework,  a leading CAV derives and implements a control strategy that forces the following HDVs to form a platoon. Using a predefined car- following model, we provided a complete, analytical solution of the CAV strategy intended for the platoon formation. We also provided a detailed analysis along with appropriate conditions under which a feasible platoon formation time exists.  Ongoing research considers the notion of optimality to derive energy-  and time-optimal platoon formation framework
by relaxing the assumption of the steady-state traffic flow. Future research should extend the proposed framework to make it agnostic of the car-following models.
%
%
%
%

\newpage
\cfoot{\thepage} 
\setlength{\textheight}{605pt}
\bibliographystyle{IEEEtran}
\bibliography{reference1,reference2,reference3,reference4,reference5,reference6,reference7,reference8,references9,references10,references11,references12,references13,references14,referenceTAC,udssc,referenceIDS,referenceRAS,platoon}

\end{document}